\documentclass[preprint,12pt]{elsarticle}

\usepackage[latin1]{inputenc}
\usepackage[T1]{fontenc}
\usepackage[english]{babel}

\usepackage{dsfont,yfonts,amsfonts,pifont}
\usepackage{amsmath,amssymb}                                                               
\usepackage{amsthm}
\usepackage{fancyhdr}
\usepackage{titletoc}
\usepackage{amsthm}

\usepackage[colorlinks=true,linktocpage=true,citecolor=blue,urlcolor=blue,bookmarksopen=true,bookmarksnumbered=true]{hyperref} 

\newtheoremstyle{theoreme}
{1cm}
{1cm}
{\em}   
{}   
{\scshape\bfseries}
{.}  
{.5em}
{}

\theoremstyle{theoreme}
\newtheorem{thm}{Theorem}
\newtheorem*{thm*}{Theorem}
\newtheorem{lem}{Lemma}
\newtheorem{prop}{Proposition}

\DeclareMathOperator{\quadro}{\scriptscriptstyle \square}
\newcommand{\transp}[1]%
{{\vphantom{#1}}^{\mathit t}{\hspace{-.05cm}#1}}

\def\D{\mathcal{D}}
\def\S{\mathcal{S}}
\def\Str{\textnormal{Str}}
\def\Tr{\textnormal{Tr}}
\def\GL{\textnormal{GL}}
\def\Id{\textnormal{Id}}
\def\Aut{\textnormal{Aut}}
\def\End{\textnormal{End}}

\journal{}

\begin{document}

\begin{frontmatter}


\ead{deolivei@iecn.u-nancy.fr}

\title{On the automorphism group of tube type real symmetric domains}


\author{Fernando De Oliveira}

\address{Institut \'Elie Cartan de Nancy (IECN), Nancy-Université, CNRS, INRIA, Boulevard des Aiguillettes, B.P. 239, 54506 Vandoeuvre-lès-Nancy, France}

\begin{abstract}
The aim of this note is to explain a generalization to the real case of a well known result on the automorphism group of an unbounded tube type symmetric domain in a complex vector space of finite dimension.
\end{abstract}


\end{frontmatter}




\section{Introduction}
\label{intro}

Let $\D=G/K$ be a (complex) \emph{bounded symmetric domain of tube type} in a finite dimension complex vector space $V$. The tangent space at the origin is a \emph{positive hermitian Jordan triple system} $(V,\{\})$ and $V$ is also endowed to a structure of \emph{semisimple complex Jordan algebra}. Moreover, there exists an \emph{euclidian real form} $V^+$ of the Jordan algebra $V$ and a \emph{Cayley transform} $\gamma$ such that $\gamma(\D)=T_\Omega$, $T_\Omega$ being the \emph{tube domain} $T_\Omega=\Omega\oplus iV^+$ where $\Omega$ means the \emph{symmetric cone} of invertible squares of $V^+$ (see \cite{Loo77}). It is well known (see \cite[Theorem X.5.6]{FK94}) that the automorphism group $\gamma\circ G\circ \gamma^{-1}$ of $T_\Omega$ is generated by $H(\Omega), N^+$ and the inversion map of the Jordan algebra $V$, where 
$$
N^+=\{x\mapsto x+iv|v\in V^+\},\quad H(\Omega)=\{g\in\GL(V^+)|g\Omega=\Omega\}.
$$
This result can be generalized to any real bounded symmetric domain of tube type. This work is essentially based on \cite{FK94} and \cite{Loo77}.

\section{Algebraic framework}

Let $V=V^+\oplus V^-$ be the \emph{Cartan decomposition} of a real semisimple Jordan algebra $V$ of finite dimension with Cartan involution $\star$. That is
$$
V^\pm=\{x\in V|x^\star=\pm x\}
$$
and the trace form $\alpha:(x,y)\mapsto\Tr(L(xy))$ is positive definite on $V^+$ and negative definite on $V^-$, where the $L(x)$ are the multiplication operators on $V$. This implies that $V^+$ is a semisimple \emph{euclidian Jordan algebra}. For $x,y\in V$, define
$$
x\quadro y:=L(xy^\star)+[L(x),L(y^\star)].
$$
Then a simple calculation shows that $V$ equiped with the triple product
$$
\{x,y,z\}:=x\quadro y(z)
$$
is a \emph{positive real Jordan triple system}, that is $(V,\{\})$ is a real Jordan triple system and the trace form $\beta:(x,y)\mapsto\Tr(x\quadro y)$ is positive definite on $V$. Observe that $\alpha(x,y^\star)=\beta(x,y)$ so if $\transp{f}$ means the adjoint operator of $f\in\End(V)$ with respect to $\alpha$, then $\transp{f^\star}$ means the adjoint of $f$ with respect to $\beta$, where $f^\star$ is define as the conjugation of $f$ by the Cartan involution $\star$.

Distinguish the two structures of $V$ by noting $V'$ the Jordan triple structure of $V$.\\

Let $P,Q$ be the quadratic representations of $V, V'$ respectively, define as 
$$
P:x\mapsto (P(x)=2L(x)^2-L(x^2)),
$$
$$
Q:x\mapsto (Q(x):y\mapsto P(x)y^\star).
$$
Let now $\Str(V), \Str(V')$ be the structure groups and $\Aut(V), \Aut(V')$ the automorphism groups of $V$ and $V'$ respectively, that is to say
\begin{align*}
\Str(V) & :=\{g\in\GL(V)|P(gx)=gP(x)\transp{g}\} ,\\
\Aut(V) & :=\{g\in\GL(V)|g(xy)=(gx)(gy)\} ,\\
&\\
\Str(V') & :=\{g\in\GL(V)|Q(gx)=gQ(x)\transp{g^\star}\} ,\\
\Aut(V') & :=\{g\in\GL(V)|g\{x,y,z\}=\{gx,gy,gz\}\}.
\end{align*}
We see easily that if $g$ is an element of $\Str(V)$ then $\transp{g}=g^{-1}P(ge)$. Also, $\Aut(V)$ is the isotropy group $\Str(V)^e$ of $e$ in $\Str(V)$, and we have
$$
\Str(V)=\{g\in\GL(V)|\jmath\circ g\circ\jmath\in\GL(V)\}
$$
(see e.g. \cite[Propositions VIII.2.4 and VIII.2.5]{FK94}). Moreover, if $O(V,\beta)$ means the orthogonal group of $V$ with respect to $\beta$, then 
$$
\Aut(V')=\Str(V')\cap O(V,\beta).
$$
Furthermore, a simple computation shows that $\Str(V')=\Str(V)$.

\section{Tube type domains and automorphism groups}

A result of O. Loos says that there exists a \emph{real bounded symmetric domain of tube type} $\D$ such that the tangent space at the origin of $\D$ is $V'$ and, the domain $\D$ is equivalent to a unbounded real symmetric domain, the \emph{tube} $T_\Omega$ define as
$$
T_\Omega:=\Omega\oplus V^-,
$$
$\Omega$ being the symmetric cone of invertible squares of $V^+$ (see \cite[Theorem 10.8]{Loo77}). The equivalence is given by the \emph{partial Cayley transform} $\gamma_e$, define on
$$
\{x\in V|e-x\text{ is invertible}\}
$$
by
$$
\gamma_e(x)=(e+x)(e-x)^{-1}=-e+2(e-x)^{-1},
$$
where $e$ means the identity element of the Jordan algebra $V$. The geodesic symmetry arround $0\in\D$ is $-\Id_V$ and then, the geodesic symmetry arround $e=\gamma_e(0)\in T_\Omega$ is the inversion map $\jmath:x\mapsto x^{-1}$ of the algebra $V$, which means that 
$$
\jmath\circ\gamma_e=\gamma_e\circ(-\Id_V).
$$

The bounded domain $\D$ is homogeneous under the action of a connected semisimple Lie group $G$ without center : $\D=G/K$, $K$ being the isotropy group of $0$ in $G$ (see \cite[11.14]{Loo77}). One has that $K$ is the identity component of $\Aut(V')$ (see \cite[Lemme 2.11]{Loo77}). Now, we can define the automorphism group $L$ of the \emph{tube domain} $T_\Omega$ as
$$
L:=\gamma_e\circ G\circ\gamma_e^{-1}.
$$
The base point of $T_\Omega$ is $e$ so $T_\Omega=L/U$ with $U=\gamma_e\circ K\circ\gamma_e^{-1}$, so $U$ is the isotropy group of $e$ in $L$. 

For $v\in V^-$, let $t_v:x\mapsto x+v$ be the translation by $v$ and define $N^+:=\{t_v|v\in V^-\}$. The translation group $N^+$ appears  as an abelian subgroup of $L$ isomorphic to $V^-$. Finally, let $G(\Omega)$ be the identity component of the group
$$
\{g\in\GL(V)|g\Omega=\Omega, g^\star=g\}=\{g\in\GL(V)|g\overline{\Omega}=\overline{\Omega}, g^\star=g\}.
$$
Then $G(\Omega)$ acts naturally on $T_\Omega$ as $a+v_-\mapsto ga+gv_-$ and thus we identifie $G(\Omega)$ with a subgroup of $L$.

We shall prove the following theorem.

\begin{thm}\label{main_thm}
The group $L$ is generated by $G(\Omega), N^+$ and $\jmath$. 
\end{thm}

\section{Affine transformations of $T_\Omega$}

Recall that an element $x$ of the Jordan algebra $V$ is invertible if and only if the operator $P(x)$ is and therefore $x^{-1}=P(x)^{-1}x$. 

The group $G(\Omega)$ acts transitively on the cone $\Omega$. Indeed, let $x\in V^+$ be an invertible element. Since $P(x)^\star=P(x^\star)=P(x)$, the restriction to $V^+$ of $P(x)$ is in $\GL(V^+)$. We have also $P(x)e=x^2\in\Omega$ and so $P(x)\Omega$ is the connected component of $x^2$ in the group $(V^+)^\times$ of invertible elements of $V^+$, none other than $\Omega$. Hence
$$
\Omega=\{x^2|x\in(V^+)^\times\}=\{P(x)e|x\in(V^+)^\times\}\subset G(\Omega)e\subset\Omega.
$$
For very good and complete explanations on euclidian Jordan algebras and symmetric cones, see \cite{FK94}.

\begin{lem}\label{lemma1}
We have the decomposition
$$
L=N^+G(\Omega)U.
$$
\end{lem}

\begin{proof}
Let $g\in L$ and $x=g\cdot e:=x+y\in\Omega\oplus V^-$. Then there exists $h\in G(\Omega)$ such that $x=he$ and thus we have
$$
x=t_y\circ h(e).
$$
The transformation $g':=h^{-1}\circ t_y^{-1}\circ g$ of $T_\Omega$ satisfies $g'(e)=e$, that is to say $g'\in U$ and 
$$
g=t_y\circ h\circ g'
$$
is the desired decomposition.
\end{proof}

The next step is the characterization of the group of affine tranformations of $T_\Omega$, which is the key of the proof of the theorem \ref{main_thm}.\\

According to \cite[Proposition 2.2]{Loo77}, the group $K$ is the normalizer of $-\Id_V$ in $G$, therefore $U$ is the normalizer of $\jmath$ in $L$ :
$$
U=\{g\in L|\jmath\circ g\circ\jmath=g\}.
$$
By consequence,
$$
U\cap\GL(V)\subset\Str(V)^e=\Aut(V)
$$
so if $g\in U\cap\GL(V)$ then $\transp{g}=g^{-1}$. Let $g=\gamma_e\circ k\circ\gamma_e^{-1}\in U$ with $k\in K$. We have
\begin{align*}
g^\star=(\gamma_e\circ k\circ\gamma_e^{-1})^\star & =\gamma_e^\star\circ k^\star\circ{\gamma_e^\star}^{-1} \\
                                                  & =\gamma_e\circ\transp{k}^{-1}\circ\gamma_e^{-1} \\
                                                  & =\gamma_e\circ k\circ P(k^{-1}e)\circ\gamma_e^{-1}.
\end{align*}
By linearity of $g$, we get
$$
0=g(0)=\gamma_ek\gamma_e^{-1}(0)
$$
thus
$$
\gamma_e^{-1}(0)=k\gamma_e^{-1}(0).
$$
Since $\gamma_e^{-1}(0)=e$, we obtain $ke=e$ and $g^\star=g$. Besides, a linear element $g$ of $U$ satisfies $\transp{g^\star}^{-1}=g$ and $ge=e$, which implies
$$
U\cap\GL(V)\subset\Aut(V')^e,
$$
$\Aut(V')^e$ being the isotropy group of $e$ in $\Aut(V')$. Observe that
$$
\Aut(V')^e=\Aut(V)\cap O(V,\beta)=\{g\in\Aut(V)|g^\star=g\}.
$$
Conversely, we have $\Aut(V')^e\subset\Aut(V)$ thus any element $g$ of $\Aut(V')^e$ verifies $\gamma_e\circ g\circ\gamma_e^{-1}$, which means
$$
\Aut(V')^e=\gamma_e\circ\Aut(V')^e\circ\gamma_e^{-1}\subset U.
$$
After all, we get
$$
U\cap\GL(V)=\Aut(V')^e.
$$
Moreover, one has $\Aut(V')^e=G(\Omega)^e$, $G(\Omega)^e$ being the isotropy group of $e$ in $G(\Omega)$. Indeed, $G(\Omega)^e\subset U \cap\GL(V)=\Aut(V')^e$ and conversely, for $g\in\Aut(V')\subset\Aut(V)$ and $x\in\Omega$ ($x=y^2$ for some $y\in(V^+)^\times$), we have
$$
gx=g(y^2)=(gy)^2\in\Omega
$$
thus $\Aut(V')^e\subset G(\Omega)^e$. We just proved the following lemma :

\begin{lem}\label{lemma2}
The following equalities hold :
$$
U\cap\GL(V)=\Aut(V')^e=G(\Omega)^e.
$$
\end{lem}

We can now give the main result of this section :

\begin{prop}\label{prop1}
Let $M\in\GL(V)$ and $v\in V$. Then the affine transformation $x\mapsto Mx+v$ belongs to $L$ if and only if $M\in G(\Omega)$ and $v\in V^-$. In particular, the connected component of the identity of the affine transformations group of $T_\Omega$ is the semi-direct product $G(\Omega)\ltimes N^+$.
\end{prop}

\begin{proof}
The part \textquoteleft if\textquoteright\ of the equivalence is obvious. For the converse, suppose that $x\mapsto Mx+v$ belongs to $L$ and first show that $M\in G(\Omega)$. According to \ref{lemma1}, exists $v\in V^-$, $h\in G(\Omega)$ and $g\in U$ such that $M=t_v\circ h\circ g=h\circ t_{h^{-1}v}\circ g$. By linearity of $h^{-1}\circ M$, the map $g$ is an affine transformation which can be written $x\mapsto g'x+u$ with $g'\in\GL(V)$ and $u\in V^-$. First, we note that the equality $ge=e$ implies $g'e=e-u\in T_\Omega$ . Now, we know that we have $g=\jmath g\jmath$, thus for every $\varepsilon>0$ one has
$$
\varepsilon g'e+u=g(\varepsilon e)=\jmath g\jmath(\varepsilon e)=\varepsilon(g'e+\varepsilon u)^{-1}.
$$
By continuity of these maps, we get $u=0$. Hence $g\in\GL(V)$ and \ref{lemma2} gives $M\in G(\Omega)$. It remains to prove that $v$ is in $V^-$. The inverse map $x\mapsto M^{-1}x-M^{-1}v$ of $x\mapsto Mx+v$ is also in $L$ and for $\varepsilon>0$,
$$
M(\varepsilon e)+v\in T_\Omega .
$$
We obtain $v\in\overline{T_\Omega}=\overline{\Omega}\oplus V^-$ through continuity of $M$. Let $p_{V^+}$ be the projection operator on $V^+$. Then $p_{V^+}(v)\in\overline{\Omega}$ and also $p_{V^+}(-M^{-1}v)\in\overline{\Omega}$. But $M=M^\star$ thus $p_{V^+}(-M^{-1}v)=-M^{-1}p_{V^+}(v)$. Moreover, since $M$ stabilizes $\overline{\Omega}$, we get $-p_{V^+}(v)\in\overline{\Omega}$ and hence $p_{V^+}(v)=0$, that is $v\in V^-$.
\end{proof}

\section{Proof of the theorem}

For the sequel, we need some additional elements. Let
$$
\Sigma:=\{x\in V|x\text{ is invertible and }x^{-1}=x^\star\} ;
$$
$\Sigma$ is not empty because it contains the orbit $Ke$. Indeed, for $k\in K$, $ke$ is invertible and we have
$$
(ke)^{-1}=\transp{k}^{-1}e=k^\star e=(ke)^\star.
$$
The set $\Sigma$ corresponds in $V'$ to the set of \emph{maximal tripotents} (see \cite[11.10]{Loo77}). For $x\in\Sigma$ we define $\S_x=Kx$ ; $\S_x$ is the connected component of $x$ in $\Sigma$.

\begin{lem}\label{lemma3}
We have
$$
\gamma_e^{-1}(V^-)=\{x\in\Sigma|e-x\text{ is invertible}\}\subset\S_{-e}.
$$
In particular, $\gamma_e^{-1}(V^-)$ is a nonempty open connected set, dense in $\S_{-e}$.
\end{lem}

\begin{proof}
Let $x\in\gamma_e^{-1}(V^-)$, $x=e-2(v+e)^{-1}=(v-e)(v+e)^{-1}$ for some $v\in V^-$. Then $x$ and $e-x$ are invertible elements and 
\begin{align*}
x^{-1}=(v+e)(v-e)^{-1} & =-(v+e)(e-v)^{-1} \\
                       & =(v-e)^\star[(v+e)^{-1}]^\star=x^\star.
\end{align*}
Conversely, let $x\in\Sigma$ be such that $e-x$ is invertible. Then
\begin{align*}
\gamma_e(x)^\star & =\gamma_e(x^\star)=-e+2(e-x^\star)^{-1} \\
                  & =-e+2(e-x^{-1})^{-1}=-e-2x(e-x)^{-1} \\
                  & =-\gamma_e(x)
\end{align*}
and hence $\gamma_e(x)\in V^-$. We get
$$
\gamma_e^{-1}(V^-)=\{x\in\Sigma|\det(e-x)\neq0\}
$$
where $\det$ means the \emph{determinant} function of $V$. The map $\det$ is polynomial so is continuous and thus, $\gamma_e^{-1}(V^-)$ is an open set in $\Sigma$. The open set $\gamma_e^{-1}(V^-)$ is also connected because $V^-$ is connected and the map $\gamma_e^{-1}$ is continuous on $V^-$. Finally, $\gamma_e^{-1}(V^-)$ contains $-e$ thus it is included in $S_{-e}$.
\end{proof}

We define now a binary relation $\intercal$ on $V$ called \emph{transversality relation}, as 
$$
x\intercal y\Longleftrightarrow\det(x-y)\neq0.
$$
For $x\in\Sigma$, we define also $x_\intercal:=\{y\in\Sigma|x\intercal y\}$. From \ref{lemma3}, $e_\intercal=\gamma_e^{-1}(V^-)\subset\S_{-e}$. If $x=ke\in\S_e=Ke$ then for all $y\in\Sigma$,
$$
\det(x-y)\neq0\Longleftrightarrow\det(e-k^{-1}y)\neq0\Longleftrightarrow k^{-1}y\in\gamma_e^{-1}(V^-)
$$
and hence
$$
x_\intercal=k\gamma_e^{-1}(V^-)\subset\S_{-x}=\S_{-e}
$$
is an open connected set, dense in $\S_{-e}$. Consequently, for all $x,y$ in $\S_e$ we have $x_\intercal\cap y_\intercal\neq\emptyset$.

\begin{lem}\label{lemma4}
For all $x,y$ in $\S_e$, exists $k\in K$ satisfying $-kx\intercal e$ and $-ky\intercal e$.
\end{lem}

\begin{proof}
We deduce from the above that for all $x,y\in\S_e$, exists $z\in\S_{-e}$ such that $x\intercal z$ and $y\intercal z$. We can also choose $k\in K$ such that $kz=-e$. The elements $-kx$ and $-ky$ are thus transverse to $e$.
\end{proof}

\begin{prop}\label{prop2}
Let $g\in U$ and $h:=\gamma_e^{-1}g\gamma_e\in K$. If $h^{-1}e\intercal e$ then $g\in N^+G(\Omega)\jmath N^+$.
\end{prop}

\begin{proof}
From \ref{lemma3}, the condition $\det(e-h^{-1}e)\neq0$ ensures the existence of an element $v$ of $V^-$ verifying $h^{-1}e=\gamma_e^{-1}(v)$. Let $x\in T_\Omega$. Then
\begin{align*}
g(x)=\gamma_e(h\gamma_e^{-1}(x)) & =-e+2(e-h\gamma_e^{-1}(x))^{-1} \\
                                 & =-e+2(h(\gamma_e^{-1}(v)-\gamma_e^{-1}(x)))^{-1}.
\end{align*}
But $h$ belongs to $K$ so
\begin{align*}
(h(\gamma_e^{-1}(v)-\gamma_e^{-1}(x)))^{-1} & =\transp{h}^{-1}(\gamma_e^{-1}(v)-\gamma_e^{-1}(x))^{-1} \\
                                            & =h^\star(\gamma_e^{-1}(v)-\gamma_e^{-1}(x))^{-1}.
\end{align*}
Thus
\begin{align*}
g(x) & =-e+2h^\star(\gamma_e^{-1}(v)-\gamma_e^{-1}(x))^{-1} \\
     & =-e-h^\star((e+v)^{-1}-(e+x)^{-1})^{-1}.
\end{align*}
Using the following Hua identity (see e.g. \cite[Exercice 5, Chapitre II]{FK94}) 
$$
a^{-1}-b^{-1}=(a+P(a)(b-a)^{-1})^{-1}
$$
with $a=e+v$ and $b=e+x$, we get
\begin{align*}
g(x) & =-e-h^\star(e+v+P(e+v)(x-v)^{-1}) \\
     & =-e-h^\star(e+v)-h^\star P(e+v)(x-v)^{-1} \\
     & =-h^\star P(e+v)(\jmath\circ t_{-v}(x))-e-h^\star(e+v).
\end{align*}
The proposition \ref{prop1} gives $g\circ t_v\circ \jmath\in N^+G(\Omega)$ and hence $g\in N^+G(\Omega)\jmath N^+$.
\end{proof}

We now have all the elements to prove the theorem \ref{main_thm}, we recall :

\begin{thm*}
The group $L$ is generated by $G(\Omega), N^+$ and $\jmath$. 
\end{thm*}

\begin{proof}
Using \ref{prop2}, it suffices to prove that every element $g\in U$ such that $\det(e-h^{-1}e)=0$ is generated by $N^+, G(\Omega)$ and $\jmath$, where $h=\gamma_e^{-1}g\gamma_e\in K$. Let $g$ be such an element. From \ref{lemma4}, we can find $k\in K$ satisfying 
$$
\det(e+k^{-1}h^{-1}e)\neq0\quad\text{and}\quad\det(e+k^{-1}e)\neq0.
$$
By \ref{prop2} we get already $\widetilde{g}:=\gamma_e\circ(-k)\circ\gamma_e^{-1}=\gamma_e\circ k\circ\gamma_e^{-1}\circ\jmath$ is generated by $N^+,G(\Omega)$ and $\jmath$. Exists also $v\in V^-$ verifying $(-hk)^{-1}e=\gamma_e^{-1}(v)$. Let now $x\in T_\Omega$. One has
\begin{align*}
g(x)=\gamma_e(h\gamma_e^{-1}(x)) & =-e+2(e-h\gamma_e^{-1}(x))^{-1} \\
                                 & =-e+2(-hk\gamma_e^{-1}(v)-h\gamma_e^{-1}(x))^{-1} \\
                                 & =-e+2(-hk(\gamma_e^{-1}(v)+k^{-1}\gamma_e^{-1}(x)))^{-1} \\
                                 & =-e-2h^\star k^\star(\gamma_e^{-1}(v)+k^{-1}\gamma_e^{-1}(x))^{-1} \\
                                 & =-e-2h^\star k^\star(\gamma_e^{-1}(v)-\gamma_e^{-1}\widetilde{g}^{-1}(x))^{-1}.
\end{align*}
Repeating  the calculus of the proof of \ref{prop2}, we obtain that $g\circ\widetilde{g}$ belongs to $N^+G(\Omega)\jmath N^+$ and hence $g$ is generated by $N^+$, $G(\Omega)$ and $\jmath$.
\end{proof}

\section{Relation between $\Str(V)$ and $G(\Omega)$}

If $V^-=\{0\}$, that is if $V=V^+$ is an euclidian Jordan algebra, then $T_\Omega=\Omega$. Also the condition $g^\star=g$ for $g\in\GL(V)$ is empty. This implies that $G(\Omega)=\{g\in\GL(V^+)|g\Omega=\Omega\}$ and
$$
\{g\in\Str(V)|g^\star=g\}=\Str(V^+)=\pm G(\Omega)
$$
by \cite[Proposition VIII.2.8]{FK94}.\\

Let now $\D$ be a complex bounded symmetric domain of tube type as in the section \ref{intro}. It is the particular case where $V^-=iV^+$, that is to say when the domain $\D$ is considered as real. The condition $g^\star=g$ for $g\in\GL(V)$ is equivalent to the existence of a unique $\widetilde{g}\in\GL(V^+)$ satisfying $g(u+iv)=\widetilde{g}u+i\widetilde{g}v$ for all $u,v\in V^+$.  Then
$$
\{g\in\Str(V)|g^\star=g\}=\{g\in\GL(V)|g^\star=g \text{ and } \widetilde{g}\in\Str(V^+)\}\cong\Str(V^+)
$$
and 
$$
G(\Omega)=\{g\in\GL(V)|g^\star=g \text{ and } \widetilde{g}\in H(\Omega)\}\cong H(\Omega).
$$
Using the identity $\Str(V^+)=\pm H(\Omega)$, we have hence again
$$
\{g\in\Str(V)|g^\star=g\}=\pm G(\Omega).
$$
In fact, the last equality is always true.

\begin{lem}\label{lemma5}
An element $x$ of $V^+$ belongs to $\Omega\cup(-\Omega)$ if and only if $P(x)$ is positive definite.
\end{lem}

\begin{proof}
Let $x\in\Omega\cup(-\Omega)$ ; $x=\pm y^2$ for some $y\in(V^+)^\times$ and thus
$$
P(x)=P(\pm P(y)e)=P(y)^2
$$
is positive definite. Conversely, let $x\in V^+$ be such that $P(x)$ is positive definite. First, $P(x)^\star=P(x)$ so the restriction to $V^+$ of $P(x)$ (ever noted $P(x)$) belongs to $\GL(V^+)$. Moreover, exists a Jordan frame $\{c_j\}$ of $V^+$ and real numbers $\{\lambda_j\}$ such that $x=\sum_j\lambda_jc_j$. If $V^+=\sum_{i,j}V_{ij}^+$ means the Peirce decomposition of $V^+$ with respect to $\{c_j\}$, then $P(x)$ is $\lambda_i\lambda_j$ on $V_{ij}^+$. The positivity condition prove that $\lambda_i\lambda_j>0$ for all $i,j$. Consequently, $\lambda_i>0$ or $\lambda_i<0$ for all $i$ and hence $x\in\Omega\cup(-\Omega)$.
\end{proof}

\begin{prop}\label{prop3}
The following equality holds :
$$
\{g\in\Str(V)\,|\,g^\star=g\}=\pm G(\Omega).
$$
\end{prop}

\begin{proof}
Let $g\in\Str(V)$ be such that $g^\star=g$ and let $x\in\Omega$. Then $gx\in V^+$ and
$$
P(gx)=gP(x)\transp{g}
$$
is positive definite. From \ref{lemma5} we get $gx\in\Omega\cup(-\Omega)$. Since $\Omega$ is connected, we find $g\Omega\subset\Omega$ or $g\Omega\subset-\Omega$. It is the same for  $g^{-1}$ so $g\Omega=\Omega$ or $g\Omega=-\Omega$. Therefore,
$$
\{g\in\Str(V)|g^\star=g\}\subset\pm G(\Omega).
$$
Let now $g\in G(\Omega)$. We can find $y\in(V^+)^\times$ such that $ge=P(y)e$ and thus $P(y)^{-1}g\in G(\Omega)^e$. So we write $g=P(y)k$ with $y\in (V^+)^\times$ and $k\in G(\Omega)^e$. As $P(y)\in\Str(V)$ verifying $P(y)^\star=P(y)$ and since $G(\Omega)^e=\Aut(V)\cap O(V,\beta)\subset\{g\in\Str(V)|g^\star=g\}$, we get $G(\Omega)\subset\{g\in\Str(V)|g^\star=g\}$. Observing that $-\Id_V\in\{g\in\Str(V)|g^\star=g\}$, we obtain the result.
\end{proof}



\bibliographystyle{model1a-num-names}
\bibliography{<your-bib-database>}

\begin{thebibliography}{00}

\bibitem[1]{FK94}
J.~Faraut and A.~Kor{\'a}nyi, \emph{Analysis on symmetric cones}, Oxford
  Mathematical Monographs, Clarendon Press, Oxford, 1994.
  
\bibitem[2]{FKKLR}
J.~Faraut, S.~Kaneyuki, A.~Kor{\'a}nyi, Q.{-}k. Lu, and G.~Roos, \emph{Analysis
  and {G}eometry on {C}omplex {H}omogeneous {D}omains}, Progress in
  Mathematics, Birkh{ä}user, Boston-Basel-Berlin, 2000.

\bibitem[3]{Loo71}
O.~Loos, \emph{Jordan {T}riple {S}ystems, {$R$}-spaces, and {B}ounded
  {S}ymmetric {D}omains}, Bull. Amer. Math. Soc. \textbf{77} (1971), no.~4,
  p.558--561.

\bibitem[4]{Loo77}
\bysame, \emph{Bounded symmetric domains and {J}ordan pairs}, Mathematical
  Lectures, University of California at Irvine, 1977.

\end{thebibliography}



\newcommand{\etalchar}[1]{$^{#1}$}
\providecommand{\bysame}{\leavevmode\hbox to3em{\hrulefill}\thinspace}
\providecommand{\MR}{\relax\ifhmode\unskip\space\fi MR }
\providecommand{\MRhref}[2]{%
  \href{http://www.ams.org/mathscinet-getitem?mr=#1}{#2}
}
\providecommand{\href}[2]{#2}

\end{document}